\documentclass[12pt]{article}

\usepackage{mathptmx} 

\newcommand{\version}{March 2, 2022}

\newcommand\tprod{\textstyle\prod}
\newcommand{\lanbox}{\hfill \hbox{$\, 
\vrule height 0.25cm width 0.25cm depth 0.01cm
\,$}}

\usepackage{amsthm,amsfonts,amsmath, amscd}

\usepackage{xcolor}

\swapnumbers
\theoremstyle{plain}

\newtheorem{thm}{THEOREM}[section]
\newtheorem{lm}[thm]{LEMMA}
\newtheorem{cl}[thm]{COROLLARY}
\newtheorem{prop}[thm]{PROPOSITION}
\newtheorem{conjecture}[thm]{Conjecture}

\theoremstyle{definition}
\newtheorem{defi}[thm]{Definition}
\theoremstyle{definition}
                        
\newcommand{\upchi}{\raise1pt\hbox{$\chi$}}
\newcommand{\R}{{\mathord{\mathbb R}}}

\newcommand{\Z}{{\mathord{\mathbb Z}}}

\newcommand{\E}{{\mathord{\cal E}}}   
\newcommand{\cA}{{\mathord{\cal A}}}

\newcommand{\cI}{{\mathord{\cal I}}}
\newcommand{\cM}{{\mathord{\cal M}}}
\newcommand{\cN}{{\mathord{\cal N}}}
\newcommand{\cU}{{\mathord{\cal U}}}
\newcommand{\cV}{{\mathord{\cal V}}}

\newcommand{\ba}{{\mathbf{a}}}

\numberwithin{equation}{section}
\begin{document}


\title{On the extension of the FKG inequality to $n$ functions}
\author{\vspace{5pt} Elliott H Lieb$^1$ and
Siddhartha Sahi$^{2,}$\footnote{Corresponding author email: sahi@math.rutgers.edu} \\
\vspace{5pt}\small{$1.$ Departments of Mathematics and Physics, Jadwin
Hall,} \\
[-6pt]
\small{Princeton University,  Princeton, NJ
  08544, USA}\\
  \vspace{5pt}\small{$2.$ Department of Mathematics, Hill Center,}\\[-6pt]
\small{Rutgers University,
110 Frelinghuysen Road,
Piscataway, NJ 08854, USA}
 }
\date{\version}

\maketitle
                                                                        
        \bigskip         
        
\footnotetext{
\copyright\, 2021 by the authors. This paper may be reproduced, in its
entirety, for non-commercial   
 purposes.}

\begin{abstract}
 
The 1971 Fortuin-Kasteleyn-Ginibre (FKG) inequality for two monotone functions on a distributive lattice is well known and has seen many applications in statistical mechanics and other fields of mathematics.
In 2008 one of us (Sahi) conjectured an extended version of this inequality for all $n>2$ monotone functions on a distributive lattice. Here we prove the conjecture for two special cases: for monotone functions on the unit square in ${\mathbb R}^k$ whose upper level sets are $k$-dimensional rectangles, and, more significantly, for {\it arbitrary} monotone functions on the unit square in ${\mathbb R}^2$. The general case for ${\mathbb R}^k, k>2$ remains open.

\end{abstract} 

\leftline{\footnotesize{\qquad Mathematics subject classification numbers:  
05A20, 26D07, 60E15, 82B20}}
\leftline{\footnotesize{\qquad Key Words: monotone functions, FKG inequality, higher correlation, symmetric group}}


                     

\section{Introduction} \label{intro}
\medskip

For functions $f,g$ on a probability space $(L,\mu)$, their expectation and correlation are defined by
\begin{equation}\label{E_12}
 E_1(f)=\E(f):=\int_L f d\mu \quad \text{ and } \quad  E_2(f,g)=\E(fg)-\E(f)\E(g).
\end{equation} 
Now suppose further that $L$ is a distributive lattice\footnote{A distributive lattice is a partially ordered set, closed under join (supremum) $\vee$ and meet (infimum) $\wedge$, such that each operation distributes over the other. A key example is the power set of a set, partially ordered by inclusion.} and that the probability measure $\mu$ satisfies
\begin{equation} \label{sup-mod}
\mu(a\vee b)\mu(a\wedge b) \ge \mu(a)\mu(b). \end{equation}
In this situtation if $f,g$ are positive monotone (decreasing) functions\footnote{In this paper we use \emph{positive} as a synonym for \emph{non-negative} and \emph{monotone} for \emph{monotone decreasing}. By reversing the partial order, our results and conjectures hold equally for monotone increasing functions. We note further that the positivity requirement on functions is redundant for the second inequality but essential for the first.} on $L$, then one has
\begin{equation} \label{fkg}
 E_1(f)\geq 0  \quad \text{ and } \quad E_2(f,g) \geq 0.
 \end{equation}
The first inequality is obvious, while the second is the celebrated FKG inequality of Fortuin-Kasteleyn-Ginibre \cite{FKG} which plays an important role in several areas of mathematics/physics. We will refer to a distributive lattice $L$ with probability measure $\mu$ satisfying \eqref{sup-mod} as an \emph{FKG poset}. 

In formulating \eqref{sup-mod} we have tacitly assumed that the poset $L$ is a discrete set. However the FKG inequality also has important continuous versions, which can be proved by discrete approximation. For example, if $Q_k=[0,1]^k$ is the unit hypercube  in $\R^k$ equipped with the partial order: $x\ge y$ iff $x_i\ge y_i$ for all $i$, then the FKG inequality holds for the Lebsegue measure, and more generally for any absolutely continuous measure whose density function satisfies \eqref{sup-mod}.

In \cite{sahi}, Sahi introduced a sequence of multilinear functionals $E_n(f_1,\ldots,f_n), \; n=1,2,3\ldots$, generalizing $E_1$ and $E_2$ (see Definition \ref{def:En} below), and made the following conjecture:
\begin{conjecture}\label{conj0}(\cite{sahi}, Conjecture 5)
If $f_1,\ldots,f_n$ are positive monotone functions on an FKG poset then
\begin{equation} \label{En-Ineq}
 E_n(f_1,\ldots,f_n)\geq 0.
 \end{equation}
\end{conjecture}
The paper \cite{sahi} proves the conjecture for the lattice $\{0,1\}\times \{0,1\}$, and for a certain \emph{subclass} of positive monotone functions on the general power set lattice $\{0,1\}^k$ equipped with a product measure. Since the functionals $E_n$ satisfy the following ``branching'' property (\cite{sahi}, Theorem 6)
\begin{equation} \label{En-rest}
 E_n(f_1,\ldots,f_{n-1},1)=(n-2)E_{n-1}(f_1,\ldots,f_{n-1}), 
 \end{equation}
  the inequalities \eqref{En-Ineq} form a \emph{hierachy} in the following sense: if $C_n$ denotes the $n$-function positivity conjecture, then $C_n$ implies $C_{n-1}$ for $n>2$.

Sahi's work was inspired by that of Richards \cite{rich}, who first had the idea of generalizing the FKG inequality to more than two functions. A natural first candidate for such an inequality is the cumulant (Ursell function) $\kappa_n$, but an easy example shows that the inequality already fails for $\kappa_3$. Nevertheless Richards \cite[Conjecture 2.5]{rich} conjectured the \emph{existence} of such a hirearchy of inequalities, although without an explicit formula for $E_n$. 

Indeed for $n=3,4,5,$ Sahi's functional $E_n$ coincides with the ``conjugate'' cumulant $\kappa_n'$ introduced by Richards \cite[ formula (2.2)]{rich}, although for $n\ge 6$ one has $E_n\ne\kappa_n'$. We note also that \cite{rich} contains two ``proofs'' of the positivity of $\kappa'_3,\kappa'_4,\kappa'_5$ -- one for a discrete lattice, and the other for a continuous analog. However it seems to us that both proofs have essential gaps. Thus, beyond the special cases treated in \cite{sahi}, Conjecture \ref{conj0} remains a conjecture, even for $n=3,4,5$. 

In the present paper we provide further evidence in support of Conjecture \ref{conj0}. We consider the continuous case of the Lebesgue measure on the unit hypercube $Q_k=[0,1]^k$ in $\R^k$, and we prove the inequalities \eqref{En-Ineq} for two additional cases:
\begin{itemize} 
\item for arbitrary positive monotone functions on the unit square in $\R^2$
\item for monotone characteristic functions of $k$-dimensional rectangles in $[0,1]^k$, and, by multinearilty of $E_n$, for functions whose level sets are (not necessarily homothetic) rectangles.
\end{itemize}
  
We treat first the case of three functions on $\R^2$ in \S \ref{3f}. This introduces several key ideas, including a reduction to a non-linear inequality involving decreasing sequences. In \S \ref{nf} we define $E_n$ for arbitrary $n$, and prove Conjecture \ref{conj0}, first for characteristic functions of $k$-dimensional rectangles, then for general monotone functions on $\R^2$; that is, we extend \S \ref{3f} to all $n >3$. This requires additional ideas involving the symmetric group $S_n$, and an intricate induction on $n$. The first two subsections of \S \ref{nf} are written in complete generality, and we hope these ideas will help in the eventual resolution of Conjecture \ref{conj0}.

Since the FKG inequality has many applications in probability, combinatorics, statistics, and physics, it reasonable to suppose that the generalized inequality will likewise prove to be useful in one or more of these areas. Although we do not have a compelling application in mind, we feel that it is important to find such an application. Indeed the right application might provide additional insight into Conjecture \ref{conj0} and perhaps even suggest a line of attack.

To end this introduction we tantalize the reader with an interesting reformulation of the inequalities $E_n \geq  0$ in terms of formal power series from \cite{sahi}. First, if $F(x)$ is a positive function on a probability space $L$, then it is natural to define the geometric mean of $F$ by the formula
\begin{equation}
  G\left( F\right) =  \exp \left( \E ( \log  F )\right).  \label{=1}
\end{equation}%
Now suppose $F(x,t)$ is a power series of the form
\begin{equation}
F(x,t)=1-{f_{1}}\left( x\right) t-{f_{2}}\left( x\right) t^{2}-\cdots .
\label{=2}
\end{equation}
Then $\log \left( F\left( x, t \right) \right)$ is a well defined power series, and formula (\ref{=1}) gives
\begin{equation} \label{=3}
G\left( F\right)  = \exp \left( \E (\log ( F(x,t) ) \right)
=1-c_{1}t-c_{2}t^{2}- \cdots 
\end{equation}%
where the constants $c_{j}$ are certain algebraic expressions in various $\E(f_{i_1}f_{i_2}\cdots f_{i_p})$.

\begin{conjecture}\label{conj} (\cite{sahi}, Conjecture 4)
If the $f_1(x), f_2(x),\ldots$ is a sequence of positive monotone functions on an FKG poset then $c_n\ge0$ for all $n$.
\end{conjecture}

It turns out that Conjectures \ref{conj0} and \ref{conj} are \emph{equivalent}. One implication has already been established in \cite[\S 3]{sahi}, and we prove the other direction in the appendix to this paper. We also refer the reader to \cite{sahi2,sahi3} for related inequalities in an algebraic setting.




\section{The inequality for three functions}   \label{3f}   


For three functions, the multilinear functional $E_n$ introduced in \cite{sahi} is given by the formula
\begin{equation}\label{E3}
E_3(f,g,h)= 2\E(fgh) + \E(f)\E(g)\E(h) -\E(f)\E(gh)-\E(g)\E(fh)-\E(h)\E(fg).
\end{equation}
We note that $E_3$ is \emph{different} from the cumulant (Ursell function) which is given by
\begin{equation}\label{kappa3}
\kappa_3(f,g,h)= \E(fgh) + 2\E(f)\E(g)\E(h) -\E(f)\E(gh)-\E(g)\E(fh)-\E(h)\E(fg).
\end{equation}

We will consider the functional $E_3$ for functions on the unit hypercube
\begin{equation}
Q_k= [0,1]^k=\{x = (x_1,\ldots, x_k)|  \  0\leq x_i\leq 1 \},
\end{equation}
 equipped with the Lebesgue measure and the usual partial order: $x\leq x'$ iff $x_i\leq x_i'$ for all $i$. We say that a real valued function $f$ on $Q_k$ is {\bf monotone} (decreasing) if $x\leq x'$ implies $f(x) \geq f(x')$. We note that the FKG inequality is usually stated for monotonically increasing functions, but this is a somewhat arbitrary choice. Indeed FKG and our theorems for decreasing functions are equivalent to the corresponding results for increasing functions. For a general FKG poset this follows by reversing the partial order, and for $Q_k$ by the change of variables $x_i\mapsto 1-x_i$. We also note that monotonicity for $Q_1$ has the usual 1-variable meaning of a decreasing function. 


 \begin{thm}\label{thm1}
If $f,g,h$ are positive monotone functions on $[0,1]^2$ then $E_3(f,g,h)\ge0$. 
\end{thm}
The generalization of Theorem \ref{thm1} to $n$ functions is given in Theorem \ref{thm3} below.  

We now reduce Theorem \ref{thm1} to characteristic functions $\chi_S,\ S\subset Q_k$.  These are defined by $\chi_S(x)=1$ if $x\in S$ and $\chi_S(x)= 0$ if $x\notin S$. We will say $S$ is {\bf monotone} if $\chi_S$ is monotone.

\begin{lm} \label{lm-char} It suffices to prove Theorem \ref{thm1} for $\chi_S,\chi_T,\chi_U$, for all monotone $S,T,U$.
\end{lm} 

{\bf Proof: } Any positive $f$ can be written as an integral
over the characteristic functions of its upper level sets. Thus, $f(x) =\int_0^\infty \xi_s (x) {\rm d}s$, with $\xi_s(x) =1 $ if $f(x)>s$ and $0$ otherwise. (See
\cite[`layer cake principle']{anal}.)  If $f$ is monotone, then $\xi_s$ is monotone for every $s$. Since $E_3$ is multi-linear in $f,g,h$, this reduces Theorem \ref{thm1} to the case of monotone characteristic functions. \lanbox

We now describe a further reduction of Theorem \ref{thm1} to a discrete family of characteristic functions. Let $\cA=\cA(m)$ be the set of decreasing $m$-tuples of integers, each between $0$ and $m$
\begin{equation}\label{Am}
 \cA(m) :=   \{a\in \Z^m \ |\ m\geq a_1 \geq \cdots \geq a_m \geq 0\}.
 \end{equation}
For each $a\in\cA$ we define a monotone subset $S_a$ of $Q_2=[0,1]^2$ as follows. Divide $Q_2$ uniformly into $m^2$ little squares, write $D_{i,j}$ for the square with top right vertex $(i/m,j/m)$, and set 
\begin{equation}\label{Sa} 
S_a =\bigcup\nolimits_{j\le a_i} D_{i,j}, \quad \chi_a=\chi_{S_a}.
\end{equation}
Then $S_a$ is a monotone subset of $Q_2$, and conversely \emph{any} monotone union of $D_{i,j}$ is of this form. 

\begin{lm} \label{lm-chara} It suffices to prove Theorem \ref{thm1} for $\chi_a,\chi_b,\chi_c$; $a,b,c\in \cA(m)$;  for all $m$.
\end{lm} 
{\bf Proof: } By Lemma \ref{lm-char} it suffices to consider monotone characteristic functions $\chi_S,\chi_T,\chi_U$. Divide $Q_2$ uniformly into $m^2$ little squares $D_{i,j}$ as before, and let $S^m,T^m,U^m$ be the unions of the $D_{i,j}$ contained in $S,T,U$, respectively; then these are monotone subsets of $Q_2$ of the form (\ref{Sa}). Moreover $\chi_{S^m},\,\chi_{S^m}\chi_{T^m}$ etc., converge to $\chi_S,\,\chi_S\chi_T$ etc., in $L^1$ as $m\to\infty$. Thus if $E_3(\chi_{S^m},\chi_{T^m},\chi_{U^m})\ge 0$ then we get $E_3(\chi_S,\chi_T,\chi_U)=\lim_{m\to\infty} E_3(\chi_{S^m},\chi_{T^m},\chi_{U^m})\ge 0.$ \lanbox


    
 \subsection{Proof of the three function inequality in two dimensions} \label{2DE3}
 
We now prove Theorem \ref{thm1} for $\chi_a,\chi_b,\chi_c$, which suffices by Lemma \ref{lm-chara}. To simplify notation, we work directly with $a,b,c$, and we define the product $ab$, expectation $\E(a)$, {\it etc.}, as follows:
 \begin{align}\label{E_i}
 (ab)_i &= \min\{a_i,b_i\},\\
 E_1(a) &= \E (a) = (a_1+ \cdots +a_m)/m^2,\\
 E_2(a,b) &= \E(ab) - \E(a) \E(b),\\
\label{E3abc} E_3(a,b,c) &= 2\E(abc) + \E(a)\E(b)\E(c)-\E(a)\E(bc)-\E(b)\E(ac)-\E(c)\E(ab).
\end{align}
Then we have $\chi_{ab}=\chi_a\chi_b$, $\E(a)= \E(\chi_a)$, $ E_2(a,b)=E_2(\chi_a,\chi_b)$, $E_3(a,b,c)=E_3(\chi_a,\chi_b,\chi_c)$. 

In particular, by the FKG inequality we get: 
\begin{lm} \label{lem2} For all $a,b$ in $\cA$ we have $E_2(a,b)\geq 0$. \hfill\lanbox\end{lm} 

To study $E_3(a,b,c)$ we consider certain perturbations of $a$. We say that $a\in \cA$ has {\bf descent} at $i$ if $a_i >a_{i+1}$, and in this case we can define three new sequences $a^{-}=a^{-,i}$, $a^{+}=a^{+,i}$, $a^\star=a^{\star,i}$, \emph{also} in $\cA$, in which the following changes, and {\it only these}, are made to $a$:
\begin{equation}\label{seven}
 a_i^- = a_{i+1},
 \quad  a_{i+1}^+ = a_i,
 ,  \quad  a_{i+1}^\star = a_{i+1}+1.
 \end{equation}
 
\begin{lm}\label{lm:apmb}
If $a$ has descent at $i$, but $b$ does not, then we have $\E(a^+b) +\E(a^-b) =2\E(ab)$.
\end{lm}
{\bf Proof: } Let $b_i=b_{i+1} = \beta$, say, then we have
\begin{align}
(a^+b)_i =(a^+b)_{i+1}&= \min\{a_i,\beta\}= (ab)_i\\  
(a^-b)_i =(a^-b)_{i+1} &=\min\{a_{i+1},\beta\} 
= (ab)_{i+1}.
\end{align} 
Since the three sequences $a^+b$, $a^-b$, and $ab$, coincide except at $i,i+1$ the result follows. \lanbox

  \begin{prop}\label{prop-pm}                  
If $a$ has descent at $i$, but $b$ and $c$ do not, then
\begin{equation}\label{eq-pm}
 E_3(a^+,b,c) + E_3(a^-,b,c) = 2 E_3(a,b,c).
 \end{equation} 
  \end{prop}
{\bf Proof: } Each term of \eqref{E3abc} has a unique factor involving $a$, which is of the form $\E(ad)$, where $d=1,b,c,bc$ is a sequence in $\cA$ that does not have descent at $i$. By Lemma \ref{lm:apmb} we get 
\begin{equation}
\E(a^+d) + \E(a^-d) = 2 \E(ad)
\end{equation}
The result now follows from formula \eqref{E3abc}.
\lanbox

\begin{lm}\label{lm:astb}
If $a,b$ have descent at $i$ and $b_{i+1}\leq a_{i+1}$, then $a^\star b =ab.$ 
\end{lm}
  {\bf Proof: } Evidently $(a^\star b)_j=(ab)_j$ for $j\ne i+1$, and since $b_{i+1}\leq a_{i+1}$ we also have  
     \begin{equation}
   (a^\star b)_{i+1} = b_{i+1} = (ab)_{i+1}.
   \end{equation} 
Thus we get $a^\star b =ab$, as claimed. \lanbox

 \begin{prop}\label{prop-st}                  
If $a$ and $b$ have descent at $i$ and $b_{i+1}\leq a_{i+1}$, then we have $a^\star b =ab$ and 
  \begin{equation}
 E_3(a^\star,b,c) \leq E_3(a,b,c) \quad \text{for all } c.
 \end{equation}
  \end{prop}
    {\bf Proof: } By Lemma \ref{lm:astb} we get $\E(a^\star b) = \E(ab), \,  \E(a^\star bc) =\E(abc)$, and it follows that
 \begin{equation} \label{eq:astbc}
  E_3(a,b,c) -E_3(a^\star,b,c) = E_2(b,c)\left[\E(a^\star) -\E(a)\right] +
            \E(b)\left[\E(a^\star c) - \E(ac) \right].
\end{equation}
Evidently we have $\E(a^\star) \geq \E(a)$ and $\E(a^\star c) \geq \E(ac)$, and by the FKG inequality we also have $ E_2(b,c) \geq 0$. Thus all terms on the right of \eqref{eq:astbc} are positive, which proves the result. \lanbox

\begin{thm}\label{thm2}
  For all $a,b,c$ in $\cA$ we have $E_3(a,b,c)\geq 0.$ 
\end{thm}
  
{\bf Proof: }  Let $\cU $ be the set of triples $(a,b,c)$ in $\cA$ for which $E_3(a,b,c)$ attains its {\it minimum}, and let $\cV$ be the subset of $\cU$ for which the quantity $\E(a)+\E(b)+\E(c)$ attains its {\it maximum}.

We claim that if $(a,b,c)\in\cV$ then $a,b,c$ are constant sequences. If this is not the case, then $a$, say, has a descent at some $i$.  If $b,c$ do not have descent at $i$ then by Proposition \ref{prop-pm} we get
$$E_3(a,b,c)=\left(E_3(a^+,b,c)+E_3(a^-,b,c)\right)/2.$$ 
By minimality, $E_3(a,b,c)\le E_3(a^\pm,b,c)$, which forces $E_3(a,b,c)=E_3(a^\pm,b,c)$. Replacing $a$ by $a^+$, we reach a contradiction since $\E(a^+) > \E(a).$

If $b$, say, also has descent at $i$, then by symmetry we may assume  $b_{i+1} \leq a_{i+1}$. Then by Proposition \ref{prop-st}, $E_3(a^\star,b,c)\le E_3(a,b,c)$, and we again reach a contradiction since $\E(a^\star)>\E(a).$ 

Now we may assume $a,b,c$ are constant sequences, and, by symmetry, further assume that
$$ a\equiv m\alpha, b\equiv m\beta, c\equiv m\gamma,\quad 0\le\alpha\le\beta\le\gamma\le 1$$ 
and it follows that \
$
E_3(a,b,c) = 2\alpha+\alpha\beta\gamma -( \alpha\beta+\alpha\beta+\alpha\gamma)
=\alpha(1-\beta)(2-\gamma)\ge 0. 
$ \lanbox

\medskip
This proves Theorem \ref{thm1} for $\chi_a,\chi_b,\chi_c$ and thus by Lemma \ref{lm-chara}, in general.




\section{The inequality for $n$ functions} \label{nf}

\subsection{The definition of $E_n$}
In this subsection and the next we work with arbitrary functions on a probability space.
We start by recalling the definition of the multilinear functional $E_n(f_1,\ldots,f_n)$from \cite{sahi}. This involves the decomposition of a permutation $\sigma$ in the symmetric group $S_n$ as a product of disjoint cycles:
\begin{equation}\label{thirteen}
\sigma = (i_1, \dots, i_p)(j_1,\dots ,  j_q) \cdots .
\end{equation}
For $\sigma$ as in \eqref{thirteen} we write $C_\sigma$ for the number of cycles in $\sigma$ and we set
\begin{equation}\label{E_sig}
E_\sigma (f^1,\ldots, f^n) = \E( f^{i_1}\cdots f^{i_p})\E(f^{j_1} \cdots f^{j_q} ) \cdots .
\end{equation}
Then the following definition is due to Sahi \cite{sahi}. 
\begin{defi} \label{def:En}
For functions  $f^1,\ldots,f^n$ on a probability space $X$ we define
\begin{equation}
\label{E_n} E_n (f^1,\ldots,f^n) =  \sum\nolimits_{\sigma \in S_n} (-1)^{C_\sigma-1}  E_\sigma (f^1,\ldots,f^n).
\end{equation}
\end{defi}
Using \eqref{E_n} one can easily verify that $E_1,E_2, E_3$ coincide with their earlier definitions. We note that the factor of $2$ in the term $2E(f^1f^2f^3)$ in formula \eqref{E3} comes from the two $3$-cycles $(123)$ and $(213)$. More generally $E_n$ will have repeated terms because $E_\sigma$ is unchanged if we rearrange the indices within a cycle. For example, for $n=4$ we have 
$$
\begin{aligned}
E_4(f^1,f^2,f^3,f^4) &=  6\E(f^1f^2f^3f^4) -2\left[ \E(f^1)\E(f^2f^3f^4) + \E(f^2)\E(f^1f^3f^4)+ \cdots \right]\\ 
& +\left[  \E(f^1)\E(f^2)\E(f^3f^4)+\E(f^1)\E(f^3)\E(f^2f^4)+\cdots \right]\\
&-\left[\E(f^1f^2)\E(f^3f^4)+\E(f^1f^3)\E(f^2f^4)+\cdots\right]-\E(f^1)\E(f^2)\E(f^3)\E(f^4).
\end{aligned}
$$
We now give an explicit formula for $E_n$ in a special case.

\begin{lm} \label{1d} Let $X=[0,1]$ be the unit interval equipped with Lebesgue measure, and let $f^i$ be the characteristic function $\chi_{[0,a_i]},\ 0\le a_i\le 1$, with 
$0\le a_1\le \cdots \le a_n\le 1$. Then we have
$$ E_n (f^1,\ldots,f^n) =a_1(1-a_2)\cdots (n-1-a_n).$$
\end{lm}

We note that the above formula implies that $E_n$ is positive, i.e. that Conjecture \ref{conj0} holds for the Lebesgue measure on $[0,1]$. While it is easy enough to give a direct proof the lemma, we prefer to postpone the proof to the next subsection where we will derive it as a consequence of a more general result. 

\subsection{Algebraic properties of $E_n$}
 
 We first prove a recursive formula relating $E_n$ to $E_{n-1}.$
 \begin{prop}\label{prop2} We have  $E_n(f^1,\ldots, f^{n-1},f)=e_1+\cdots+e_{n-1} -e_n$ where 
 \begin{equation}\label{seventeen}
e_i = 
 \begin{cases} 
 E_{n-1}(f^1,\dots,f^i f, \dots, f^{n-1})  &\text{if } 1\le i\le n-1, \\
 E_{n-1}(f^1,\dots, f^{n-1})\, \E(f)         &\text {if } i=n.
 \end{cases}
\end{equation}
\end{prop}

{\bf Proof:} We write $f=f^n$ and consider the expression (\ref{E_n}) for $E_n(f^1,\dots, f^n)$ as a sum over the symmetric group $S_n$. We decompose $S_n$ as a disjoint union
\begin{equation}
S_n=S^{(1)}\cup \cdots\cup S^{(n)}, \quad S^{(i)}=\{\sigma\in S_n\mid \sigma(i)=n\}.
\end{equation}
Then $S^{(n)}$ is a subgroup of $S_n$, naturally isomorphic to $S_{n-1}$. By (\ref{E_n}) we have
\begin{equation}
E_n (f^1,\dots, f^n)=\Sigma^1+\cdots+\Sigma^n,\quad \Sigma^i=\sum\nolimits_{\sigma \in S^{(i)}} (-1)^{C_\sigma-1}  E_\sigma (f^1,\ldots,f^n).
\end{equation}

To study the $\Sigma^i$ we consider the map $\sigma \mapsto \overline{\sigma}$ defined by dropping $n$ from the cycle decomposition of $\sigma$. Thus for $n =5$ we have $(13)(245) \mapsto (13)(24)$, \ $(12)(34)(5)\mapsto (12)(34)$, \ {\it etc.}
Then $\sigma\mapsto\overline{\sigma}$ defines a bijection from \emph{each }$S^{(i)}$ to $S_{n-1}$. If $\sigma$ is in $S^{(i)}$ and $i\ne n$ then $i$ and $n$ occur in the same cycle of $\sigma$, and dropping $n$ does not change the cycle count. This gives
$$ C_\sigma= C_{\overline{\sigma}},\quad E_\sigma(f^1,\dots, f^{n-1},f) = E_{\overline{\sigma}}(f^1,\dots,f^i f, \dots, f^{n-1}), $$ which implies $\Sigma^i=e_i$. If $\sigma$ is in $S^{(n)}$ then $(n)$ occurs as a separate cycle in $\sigma$ and we get 
$$C_\sigma= C_{\overline{\sigma}}+1,\quad E_\sigma(f^1,\dots, f^{n-1},f) = E_{\overline{\sigma}}(f^1,\dots, f^{n-1})\E(f),$$
which gives  $\Sigma^n=-e_n$. This proves the Proposition. \lanbox

Lemma \ref{1d} is now an easy consequence.

{\bf Proof of Lemma \ref{1d}} Let $f_i=\chi_{[0,a_i]}$. Since $a_i\le a_n$ for all $i$ we get 
$$f_if_n=\chi_{[0,a_i]}\chi_{[0,a_n]}=\chi_{[0,a_n]}=f_i.$$
Now applying Proposition \ref{prop2} with $f=f_n$ we deduce that 
$$ E_n (f^1,\ldots,f^n) =\left((n-1) -\E(f^n)\right) E_{n-1} (f^1,\ldots,f^{n-1}) =(n-1-a_n)E_{n-1} (f^1,\ldots,f^{n-1}).$$
The result follows by a straightforward induction on $n$. \lanbox

We next establish a useful formula for the partial sum $P_c$ of $E_n$ over the set of permutations containing a fixed cycle $c$.
\begin{prop} \label {prop-Pc} Let $S^c$ denote the set of permutations $\sigma \in S_n$ that contain a fixed cycle $c=(i_1,\ldots,i_p)$ and let $J_c= \{j_1,j_2,\ldots\} =\{1,\ldots,n\} \setminus\{i_1,\ldots,i_p\}$, then we have
\begin{equation}\label{=Pc}
P_c:=\sum_{\sigma \in S^c} (-1)^{C_\sigma-1}E_\sigma(f^1,\ldots, f^n)= \begin{cases} -\E(f^{i_1}\cdots f^{i_p}) E_{n-p}(f^{j_1}, f^{j_2},\ldots) &\text{ if } p<n,\\ \qquad\qquad \E(f^1\cdots f^n) &\text{ if } p=n.
\end{cases}
\end{equation}
\end{prop}
{\bf Proof:} The set $S^c$ consists of a single permutation if $p=n$. Otherwise it consists of permutations of the form $\sigma =c \cdot \tau$ where $\tau$ is a permutation of $J^c$. Evidently the number of cycles in $\sigma$ and $\tau$ are related by $C_\tau =C_\sigma-1$. Thus in this case we have
\begin{equation} \label{Est}
(-1)^{C_\sigma -1}E_\sigma (f^1\cdots f^n)=-\E(f^{i_1}\cdots f^{i_p}) (-1)^{C_\tau -1} E_\tau(f^{j_1}, f^{j_2},\ldots).
\end{equation} 
Now the result follows by summing \eqref{Est} over $\tau$.
\lanbox

\medskip

\subsection{Proof of the $n$ function inequality for rectangles in any dimension}\label{nRk}

By a rectangle in dimension $k$, or a $k$-rectangle, we mean a subset of $[0,1]^k$ of the form
$$ [0,r_1]\times\cdots\times [0,r_k], \quad 0\le r_1,\ldots,r_k\le 1.$$

\begin{thm} If $f^i$ are characteristic functions of $k$-rectangles then $E_n(f^1,\ldots,f^n)\ge 0.$
\end{thm}
{\bf Proof:} We proceed by induction on $k\ge1$, and for a given $k$ by induction on $n\ge1$. The base cases $k=1$ and $n=1$ are straightforward, the former by Lemma \ref{1d}. Thus we may assume $k>1$ and $n>1$, and we can write 
$$f^i =g^i\times \chi_{[0,a_i]}$$ 
where $g^i$ is the characteristic function of a $(k-1)$-rectangle. By symmetry of $E_n$ we may assume \begin{equation}\label{=ai}
0\le a_1 \le \cdots \le a_n\le 1.
\end{equation} 

We note that the assumption \eqref{=ai} on the $a_i$ means that we have 
\begin{equation}\label{=Ec}
\E(f^{i_1}\cdots f^{i_p})=a_l \E(g^{i_1}\cdots g^{i_p}),\quad l=\min\{i_1,\ldots,i_p\}.
\end{equation}
Moreover it follows from \eqref{E_n} and \eqref{=Ec} that if $a_2=\cdots=a_n=1$ then we have  
\begin{equation} \label{=fg}
E_n(f^1,\ldots,f^n)=a_1 E_n(g^1,\cdots,g^n).
\end{equation}

We now fix an index $i>1$ and let $C(i)$ denote all set of all cycles containing $i$ then we have 
$$ E_n(f^1,\ldots,f^n) = \sum\nolimits_{c\in C(i)} P_c$$
where $P_c$ is as in Proposition \ref{prop-Pc}. If $i$ is not minimal in $c$ then $P_c$ is independent of $a_i$ by \eqref{=Ec}.  If $i$ is minimal in $c$ then $1\not\in c$; hence $c$ has length $p<n$ and by \eqref{=Pc} and \eqref{=Ec} we get 
$$ P_c=-a_i b_c, \quad b_c=\E(g^{i_1}\cdots g^{i_p})E_{n-p}(f^{j_1}, f^{j_2},\ldots).$$
By induction on $n$ we have $b_c\ge 0$ for such $c$.  This means that $E_n(f^1,\ldots,f^n)$ decreases as we increase $a_2\ldots,a_n$ subject, of course, to condition \eqref{=ai}. In particular, $E_n$ decreases as we successively increase 
$$ 
a_n\nearrow1, \quad a_{n-1}\nearrow1,\quad \ldots, \quad a_2\nearrow1.
$$
By \eqref{=fg} we get $E_n(f^1,\ldots,f^n)\ \ge \ a_1 E_n(g_1,\cdots,g_n)$, which is positive by induction on $k$. \lanbox
\medskip

If $f$ is the characteristic function of a rectangle, then any level set of $f$ is either the same rectangle, or empty. However, using the layer-cake principle \cite{anal} and multilinearlity as in the proof of Lemma \ref{lm-char} we obtain the following immediate extension of the previous result.  
\medskip
\begin{cl} If $f^1,....,f^n$ are positive, monotone functions whose
 level sets are (not necessarily homothetic) rectangles then $E_n(f^1,.....,f^n) \geq 0$.
\lanbox  \end{cl}


\subsection{Proof of the $n$ function inequality in two dimensions}

 Our main result is as follows.

\begin{thm}\label{thm3}
If $f^1,\ldots,f^n$ are positive and monotone on $[0,1]^2$ then $E_n (f^1,\ldots,f^n) \ge0$. 
\end{thm}

As before we can deduce this from the special case of $\chi_a$ as in (\ref{Sa}).

\begin{lm} \label{lm-charn} It suffices to prove Theorem \ref{thm3} for $\chi_{a^1},\ldots,\chi_{a^n}$, $a^i\in \cA(m)$, for all $m$.
\end{lm} 
{\bf Proof:} This is proved along the same lines as Lemmas \ref{lm-char} and \ref{lm-chara}. \lanbox

In this section we work with $\cA=\cA(m)$ and to simplify notation, for $a^1, \dots, a^n$ in $\cA$, we set
\begin{align}\label{a_sig}
E_\sigma (a^1,\ldots, a^n) &= \E(\chi_{a^{i_1}},\ldots,\chi_{a^{i_p}})\E(\chi_{a^{j_1}},\ldots,\chi_{a^{j_q}} ) \cdots ,\\
\label{a_n} E_n (a^1,\ldots,a^n) &=  \sum\nolimits_{\sigma \in S_n} (-1)^{C_\sigma-1}  E_\sigma (a^1,\ldots,a^n).
\end{align}
Then we have $E_n(\chi_{a^1},\ldots,\chi_{a^n})=E_n (a^1,\ldots,a^n).$
        
To study the positivity of $E_n$, we first consider a special case.
\begin{prop}\label{const} If $a^i\equiv m\alpha_i$ are constant sequences, with $ 1\le \alpha_1\le\cdots\le\alpha_n\le 0$, then
\begin{equation} \label{E-alpha}
 E_n(a^1,\ldots,a^n)=\alpha_1(1-\alpha_2)\cdots (n-1-\alpha_n).
 \end{equation} 
\end{prop}

{\bf Proof:} Let $L_n=E_n(a^1,\ldots,a^n)$. Since $\alpha_i\le\alpha_n$ we have $a^ia^n=a^i$ for all $i$. Thus we get
$$L_n= \sum\nolimits_{i=1}^{n-1} L_{n-1}-L_{n-1}\E(a^n)=(n-1-\alpha_n)L_{n-1},$$
by Proposition \ref{prop2}. Now \eqref{E-alpha} follows by induction on $n$, the case $n=1$ being obvious. \lanbox
 
We now prove the generalization of Proposition \ref{prop-pm}. 

\begin{prop}\label{prop3} 
If $a$ has descent at $i$, but $a^1,\dots,a^{n-1}$ do not, then we have
\begin{equation}\label{eq:a-pm-n}
2E_n(a^1,\dots,a^{n-1},a)= E_n(a^1,\dots,a^{n-1},a^+) +  E_n(a^1,\dots,a^{n-1},a^-).
\end{equation}
\end{prop}
\smallskip

{\bf Proof:} This is proved for each term $E_{\sigma}$ in \eqref{a_n}, in exactly the same way as Proposition \ref{prop-pm}, by applying Lemma \ref{lm:apmb} to the unique factor of $E_{\sigma}$ involving $a=a^n$ in \eqref{a_sig}.
\lanbox

\medskip
We shall prove the next three theorems \emph{together} by induction on $n$.

\begin{thm}\label{thm3a} 
 If $a^1,\dots, a^{n-2},b$ are in $\cA$; $S$ is a subset of $Q_2$; and $\chi_b\chi_S =0$, then 
 \begin{equation} \label{An}
E_n(\chi_{a^1},\dots, \chi_{a^{n-2}},\chi_b,\chi_S)\le0.
 \end{equation}
\end{thm}

\begin{thm}\label{thm3b}
If $a^1, \dots, a^{n-2},b,c$ are in $\cA$; $b,c$ have descent at $i$; and $b_{i+1}\leq c_{i+1}$ then 
\begin{equation}\label{Bn}
E_n(a^1, \dots, a^{n-2},b,c^\star) \leq E_n(a^1, \dots, a^{n-2},b,c).
 \end{equation}
\end{thm}

\begin{thm}\label{thm3d}
For all $a^1,\ldots, a^n$ in $\cA$ we have 
\begin{equation} 
\label{Cn}
E_n (a^1,\ldots, a^n) \geq0.
\end{equation}
\end{thm}

{\bf Proof:} Let us write $A(n)$, $B(n)$ and $C(n)$ for the assertions of Theorems \ref{thm3a}, \ref{thm3b} and \ref{thm3d}. Then $A(1),\ B(1)$ are vacuously true, while $C(1)$ is evident. Therefore it suffices to prove the implications $A(n-1)\wedge C(n-1)\implies A(n)$, and $A(n) \implies B(n) \implies C(n)$, for all $n\ge2$.

\medskip
$A(n-1)\wedge C(n-1)\implies A(n)$: \ By assumption we have $\chi_b\chi_S=0$, also we have $\chi_{a^i}\chi_S=\chi_{S^i}$ where $S^i=S\cap S_{a^i}$. Thus by Proposition \ref{prop2} we get
\begin{align*}
E_n(\chi_{a^1}&,\dots, \chi_{a^{n-2}},\chi_b,\chi_S)=e_1+\cdots+e_{n-2}+e_{n-1}-e_n,\\
\text{where }\  &e_i := E_{n-1}(\chi_{a^1},\dots, \chi_{S^i},\ldots,\chi_{a^{n-2}}, \chi_b),\quad i\le n-2,\\
&e_{n-1} := E_{n-1}(\chi_{a^1},\dots, \chi_{a^{n-2}},0)\\
&e_n := E_{n-1}(\chi_{a^1},\dots, \chi_{a^{n-2}},\chi_b) \E(\chi_S).
\end{align*}
Now $e_n\ge 0$ by $C(n-1)$, and $e_{n-1}=0$ by \eqref{a_sig} and \eqref{a_n}. Also  $\chi_b\chi_{S^i}= (\chi_b\chi_S)\chi_{a^i}=0$, and so by symmetry we can apply $A(n-1)$ to conclude $e_i\le0$ for $i\le n-2$. This implies $A(n)$,\ \eqref{An}.

\medskip

$A(n)\implies B(n)$:  Define $S_c, S_{c^\star}$ as in \eqref{Sa} and put $S=S_{c^\star}\setminus S_c$, then by Lemma \ref{lm:astb} we have 
$$\chi_S\chi_b=(\chi_{c^\star}-\chi_c)\chi_b =\chi_{c^\star b}-\chi_{cb}=0.$$ 
Thus by $A(n)$, \eqref{An}, we get $E_n(\chi_{a^1},\dots, \chi_{a^{n-2}},\chi_b,\chi_{c^\star}-\chi_{c})\le0$ which implies $B(n)$, \eqref{Bn}.

\medskip
$B(n) \implies C(n)$: This argument is similar to the proof of Theorem \ref{thm2}. Let $\cM$ be the set of $n$-tuples $\ba= (a^1,\dots, a^n)$ in $\cA$ for which $E_n(\ba)$ achieves its \emph{minimum}, and let $\cN$ be the subset of 
$\cM$ for which $\lambda(\ba)=\E(a^1)+\dots + \E(a^n)$
achieves its \emph{maximum} on $\cM$. We claim that for $\ba$ in $\cN$ each $a^i$ is a constant sequence; by Proposition \ref{const} this clearly implies $C(n)$, $E_n(\ba)\ge0$. 

If the claim is not true then one of the sequences has a descent at some $i$. First suppose that only one sequence, by symmetry $a^n=a$, has descent at $i$. By Proposition \ref{prop3} and minimality of $E_n(\ba)$ we deduce $E_n(\ba)=E_n(a^1,\dots,a^{n-1},a^\pm)$. Thus replacing $a$ by $a^+$ preserves $E_n(\ba)$ but  increases $\lambda(\ba)$, which is a contradiction. If two sequences have descent at $i$, then by symmetry we may assume these are $a^{n-1}=b$, $a^n=c$ with $b_{i+1}\le c_{i+1}$. Now $B(n)$, \eqref{Bn}, implies that replacing $c$ by $c^\star$ does not increase $E_n(\ba)$ but it does increase $\lambda(\ba)$, which is a contradiction. \lanbox

\medskip
This proves Theorem \ref{thm3} for $\chi_{a^i}$ and thus, by Lemma \ref{lm-charn}, in general.

\section*{Acknowledgement}
This work was partially supported by NSF grants DMS-1939600, DMS-2001537, and Simons foundation grant 509766. The hospitality of the Institute for Advanced Study is gratefully acknowledged. 

\section*{Data Availabilty} Data sharing is not applicable to this article as no new data were created or analyzed in this study.


\section{Appendix}
\setcounter{section}{1}
\renewcommand{\thesection}{\Alph{section}} 

In this appendix we prove the equivalence of Conjecture \ref{conj0} and Conjecture \ref{conj}. We start by recalling some basic facts about partitions and permutations. 

A partition $\lambda$ of $n$, of length $l$, is a weakly decreasing sequence of positive integers 
$$\lambda _{1}\geq \lambda _{2}\geq \cdots \geq \lambda _{l}>0, \quad \text{ such that } \quad \lambda _{1}+\cdots +\lambda _{l}=n;$$ we say that the $\lambda_j$ are the \emph{parts }of $\lambda$, and we write $l\left( \lambda \right)=l$ and $\left\vert \lambda\right\vert=n$.

The conjugation action of $S_n$ permutes the indices in the cycle decomposition \eqref{thirteen} of an element $\sigma$. Thus the class of $\sigma$ is uniquely determined by its ``cycle type'', {\it i.e.} the partition $\lambda$ whose parts are the cycle lengths of $\sigma$, arranged in decreasing order. Moreover if $m_{i}=m_{i}\left( \lambda \right) $ denotes the number of parts of size $i$, then the conjugacy class of cycle type $\lambda$ contains $n!/z_{\lambda}$ elements where
\begin{equation}\label{z-lam}
z_{\lambda }=\tprod\nolimits_{i\geq
1}i^{m_{i}}\left( m_{i}\right) !
\end{equation}%

For a function $f$ on a probability space, we define its \emph{moments }by the formula
\begin{equation}
p_d(f)=\E(f^d) \quad \text{ and } \quad p_\lambda(f)=p_{\lambda_1}(f)\cdots p_{\lambda_l}(f).
\end{equation}
 
\begin{lm} \label{lm:E-n} We have $E_n(f,\ldots,f) =n!\sum\nolimits_{|\lambda|=n}{(-1)^{l(\lambda)-1}}z_\lambda^{-1}p_\lambda(f).$
\end{lm}

{\bf Proof:} If $\sigma$ is of class $\lambda$, then the number of disjoint cycles in $\sigma$ is $l(\lambda)$ and by \eqref{E_sig} we have $E_\sigma(f,\ldots,f) = p_\lambda(f)$. Thus the sum \eqref{E_n} for $E_n(f,\ldots,f)$ is constant over conjugacy classes, with class $\lambda$ contributing $n!/z_\lambda$ identical terms. This implies the result. \lanbox

If $f$ is as above and $u$ is a parameter then we can define the formal logarithm 
\begin{equation}\label{log-u}
\log(1-uf)= -\sum\nolimits_{i\ge 1} u^i f^i/i \ .
\end{equation}
 
 \begin{prop} \label{Euf} 
 We have $\exp\left(\E\left(\log(1-uf)\right)\right)=1-\sum_{n\ge1} u^n E_n(f,\ldots,f)/n!\ .$
 \end{prop}
{\bf Proof:} Let  $Z= \E\left(\log(1-uf)\right)$ then by \eqref{log-u} we have
\begin{equation} Z =-\sum\nolimits_{i\ge 1} u^i p_i(f)/i
\end{equation}
Writing $p_k=p_k(f)$ and $p_\lambda=p_\lambda(f)$ for simplicity, we get
\begin{equation}\label{p-id}
\exp(Z)=\prod\nolimits_{i\ge1} \sum\nolimits_{m_i\ge0} (-1)^{m_i}(u^ip_i)^{m_i}/i^{m_i}m_i! =\sum\nolimits_\lambda (-1)^{l(\lambda)}z_\lambda^{-1} p_\lambda u^{|\lambda|}.
\end{equation}
Now the result follows from Lemma \ref{lm:E-n}. \lanbox
\medskip

\begin{prop}\label{prop-p}  If $f_1,f_2,\ldots$ are functions  on a probability space then we have 
$$1-\exp\left( \E\left(\log\left(1-\sum\nolimits_if_it^i\right) \right)\right) =\sum\nolimits_{n\ge1}\sum\nolimits_{i_1,\ldots,i_n}E_n(f_{i_1},\ldots,f_{i_n})t^{i_1+\cdots+i_n}/n! \ .$$   
\end{prop}
{\bf Proof:} Let us write $A=f_1t+f_2t^2+\cdots$, then by Proposition \ref{Euf} we get
$$1-\exp\left(  \E(\log(1-A) )\right) =\sum\nolimits_{n\ge1} E_n(A,\ldots,A)/n! \ ,$$
and by multinearity of $E_n$ we have
$E_n(A,\ldots,A)=\sum\nolimits_{i_1,\ldots,i_n}E_n(f_{i_1},\ldots,f_{i_n})t^{i_1+\cdots i_n}.$
\lanbox

\begin{thm} For a set of functions $\cI$ on a probability space, the following are equivalent
\begin{enumerate}
\item For all $n$, we have $E_n(f_1,\ldots,f_n)\ge0$ if $f_1,\ldots,f_n\in\cI$. 
\item The power series $1-\exp\left( \E(\log(1-\sum_if_it^i))\right)$ has positive coefficients if $f_1,f_2,\ldots \in \cI.$
\end{enumerate}
\end{thm}
{\bf Proof:} The first statement implies the second by Proposition \ref{prop-p}. The converse was proved in \cite{sahi}, but we recall it here for completeness. Let $p_1,p_2,\ldots,p_n$ be the first $n$ primes; define 
$$k=p_1p_2\cdots p_n,\quad k_j=k/p_j,\quad N=k_1+\cdots +k_n,$$ 
and consider possible solutions of the equation $s_1k_1+\cdots+s_nk_n = N$ where $s_1,\ldots, s_n$ are integers $\ge0$. If some $s_j$ were $0$ then $p_j$ would divide the left side but not the right; thus we must have all $s_j>0$ and hence that $s_1=\cdots=s_n=1$. Now it follows from Proposition \ref{prop-p} that the coefficient of $t^N$ in the power series $1-\exp\left( \E(\log(1-\sum_{j=1}^n f_jt^{k_j}))\right)$ is precisely $E_n(f_1,\ldots,f_n)$. Thus the second statement implies the first. \lanbox
\medskip

The previous theorem proves the equivalence of Conjectures \ref{conj0} and \ref{conj}. In particular, our Theorem \ref{thm3} implies Conjecture \ref{conj} for the Lebesgue measure on the unit square in $\R^2.$

emails:\ \ lieb@princeton.edu\\
 \phantom{aaaaaaaiiaa}   sahi@math.rutgers.edu

\end{document}